\documentclass[11pt]{article}

\usepackage[margin=1in]{geometry}
\usepackage{setspace}
\usepackage[utf8]{inputenc}
\usepackage{parskip}
\usepackage{amsmath,amssymb}
\title{The Antinomy of the Liar}
\date{}
\usepackage{parskip}
\usepackage{hyperref}
\usepackage{scalerel,amssymb}

\usepackage{biblatex} %Imports biblatex package
\newtheorem{innercustomthm}{Sentence}
\newenvironment{customthm}[1]
  {\renewcommand\theinnercustomthm{#1}\innercustomthm}
  {\endinnercustomthm}
  \usepackage[utf8]{inputenc}

\begin{document}
%\maketitle
\begin{titlepage}
   \begin{center}
       \vspace*{1cm}

       \textbf{The Antinomy of the Liar}

       {Helena Jorquera Riera}\\
       %3rd year undergraduate\\
       University of Southampton\\
       June 2022\\
       
        \vspace{1cm}
         %I Helena Jorquera represent that this is my own novel work, I have recently completed an undergraduate degree with a minor in philosophy and am not enrolled in, nor have I ever been enrolled in, a graduate program in philosophy.   

   \end{center}
\end{titlepage}
%\maketitle
In this essay, I present a \textit{dialetheic} solution to the \textit{Antinomy of the Liar} and I evaluate the objection that, if the argument in \textit{Curry's paradox} is valid, accepting this solution forces us to accept an analogous solution to Curry's, which I show to be unacceptable. Thus, to accept this solution, we have to reject the argument in Curry's. I discuss some arguments defending this view. %and an objection which, although can be resisted for the Liar, cannot, be resisted for Curry's paradox, and hence, forces to attempt to solve the latter in some other way, which I will ultimately deem implausible.

Consider the sentence

\begin{customthm}{L}\label{liarL}
Sentence \ref{liarL} from \cite{thispaper} is false.
\end{customthm}
 Or, writing $F(x)$ for “$x$ is false" and $T(x)$ for “$x$ is true":
 $$\ref{liarL}=F(\ref{liarL})$$
Using the \textit{T-Schema}
\begin{equation}\label{T}\tag{T}T(x) \iff x\end{equation}
we try to determine the truth-value of \ref{liarL}: 
\begin{equation}\label{L}\end{equation}
%[F] From X is not true we can deduce that X is false and vice-versa. From X is not false we can deduce that X is true and vice-versa.
0. $\ref{liarL}= F(\ref{liarL})$ (given)\\
1.1 $T(\ref{liarL})$ (Assumption) \\
1.2 $T(F(\ref{liarL}))$ (0, 1.1; substitution of identity)\\
1.3 $F(\ref{liarL})$ (1.2; \ref{T} and modus ponens (MP henceforth))\\
1.4. $ T(\ref{liarL}) \wedge F(\ref{liarL})$ (1.1, 1.3; conjunction)\\
%1.4 $L$ is true and L is false (by 1.1 and 1.3 and conjunction)\\
2.1 $F(\ref{liarL})$ (Assumption)\\
2.2 $\ref{liarL}$ (2.1; substitution of identity)\\
2.3. $T(\ref{liarL})$ (2.2; \ref{T} and MP)\\
2.4. $ T(\ref{liarL}) \wedge F(\ref{liarL})$ (2.1, 2.3; conjunction)\\
3. $ T(\ref{liarL}) \wedge F(\ref{liarL})$ (1.1-1.4, 2.1-2.4)
%One solution is to say that L is neither true nor false. However, consider the strengthened liar:

%$L=\sim T(L)$

%for L saying that it is neither true nor fals edoes not work....
%Accepting that L is neither true nor false forces us to accept that it is also true, a contradiction, hence we reject that it is neither true nor false.

Where step 3 comes from \ref{liarL} having to be true or false. We get \ref{liarL} being both true and false, a contradiction. This is the \textit{Antinomy of the Liar}.

The \textit{dialetheic} solution to the antinomy defended by Priest \cite{priest} is to accept some contradictions to be true and (\ref{L}) as evidence that \ref{liarL} is true and false. %The argument in ... is invalid if we accept dialatheism, but we can determinr the truth value of L like so:

%assume that it is true then it is true and false, assume that it is false, then it is true and false, assume it is not true then....

%But then we discharged L twice for each assumption

%(since L must be at least one of {true, false, not true or no false} it must be that L is true and false.)
%If L is true things are as L says, L says that L is false, hence L is false, a contradiction. Hence we reject the assumption that L is true, so L is false. But since L says that L is false, things are as L says hence L is true. Thus, L is true and false.

%This is the antinomy of the liar (AL).

%One solution, defended by Priest in \cite{priest}, is to say that we should accept the argument in (\ref{L}) as proof that \ref{liarL} is both true and false, and so, accept dialetheism, the view that there are true contradictions.

%If we accept dialetheism we should reject RAA for it is possible for something to be both true and false. Hence, after all, the argument in liar is invalid. Still, the only option for the truth value of L which does not yield a contradition, (but it does), is that it is both true and false.
%Accepting dialetheism contradicts the law of non contradiction, which might be stated as 

%[LNC] $\msquare \sim(a \vee \sim a)$

%However, as well as contradicting it it might also confirm it. Since, according to dialetheism, it is possible to have x and not x, we can have (a and not a) as well as not(a and not a).

%Since L is both true and false and also not L is both true and false we have, in particular, that L is true and not L is false, hence (L and not L) is false, so LNC holds. 

One objection to dialetheism is that at least prima facie, if we have a sentence $A$ which is both true and false, or equivalently if we remain within classical logic, both $A$ and its negation $\sim A$ are true, then we can conclude any other sentence $B$ is true. This is the \textit{principle of explosion}:
\begin{equation}\label{explosion}\tag{PE}A, \sim A \vdash B\end{equation}
and it can be proven like so:
\begin{equation}\label{PE}\end{equation}
1. $A$ (given)\\
2. $A \vee B$ (1; disjunction introduction)\\
3. $\sim A$ (given)\\
4. $B$ (2, 3; disjunctive syllogism)

%where $F(A)\iff  T(\sim A)$ and $\sim A$ is the negation of $A$, hence PE follows.

However, in step 4, disjunctive syllogism is not justified  if we accept dialetheism; it might be (and it is) that A, as well as false, is true. Hence, \ref{explosion} has not been proven and it fails to prove dialetheism wrong if we reject disjunctive syllogism. 

This also shows that if we accept dialetheism, we have to reject disjunctive syllogism as well as other rules of inference.%in fact we reject much more
%less 800 words expostion up to present objection
%If A and ~A are each both true and false then A \vee B, ~A does not yield B, since A is, as well as not A, true. Neither does A being false and A or B yield B, since A is, as well as false, true. %can be sumarised if say that F(A)=  T(not A)

%"Even granting that there is an operator, say, ∗, which behaves as dialetheists claim (namely, such that in particular in some cases A is true together with ∗A), it is still perfectly possible to define a negation with all the properties of classical negation; in particular, the property of being explosive."

It can be argued that even if negation in English is correctly interpreted as a dialetheic negation, where both $A$ and $\sim A$ can be true, we can still define a connective, $-$, which cannot be interpreted as such.% and which satisfies an analogous version of \ref{explosion}. For such a connective we can restate the liar as $C=-C$, and the antinomy will be restored.

Priest \cite[chapter 4.7]{priest} gives the example of defining $-A$ as $A \rightarrow  \perp$ where $\perp$ is some unacceptable conclusion such as $\forall x T(x)$ i.e. everything is true. This is equivalent to classical negation and it satisfies an analogous version of \ref{explosion}, $A, - A \vdash B$.
%\begin{equation}\label{PEcurry}\end{equation}
%1. $A$ (given)\\
%2. $-A$ (given)\\
%2. $A \rightarrow  \perp$ (given)\\
%4. $A\rightarrow \perp$ (2, 3; subst. of id.)\\
%3. $\perp$ (1, 2; MP)

We can define a sentence $C$ to be such that

%C=-C= C $\rightarrow$ $\perp$

%C iff (T(C) $\rightarrow$ $\perp$)

%or, 

$C \iff (C \rightarrow \perp$)

which generates \textit{Curry's paradox}:
\begin{equation}\label{curry}\end{equation}
1. $C \iff (C \rightarrow$ $\perp$) (given)\\
2.1 $C$ (Assumption for conditional proof (CP henceforth))\\
2.2 $C \rightarrow$ $\perp$ (1, 2.1; MP)\\
2.3 $\perp$ (2.1, 2.2; MP)\\
3. $C \rightarrow$ $\perp$ (2.1, 2.3; CP)\\
4. $C$ (1, 3; MP)\\
5. $\perp$ (3, 4; MP)

We conclude $\perp$, which is unacceptable.

%MP mayus?
%Arguing that $C$ is both true and false is not a solution unless we accept $\perp$, so the dialetheic solution for \ref{liarL} does not work for $C$.% and if we still accept it we are forced to have dissimilar solutions for two paradoxes with a very similar structure.

%(or, equivalently, that both \ref{liarL} and $\sim$\ref{liarL} are true)
If we accept (\ref{L}) as evidence that \ref{liarL} is true and false then we should accept (\ref{curry}) as evidence for both $C$ and $-C$ being true (since we have $C$ in step 4 and $-C$ in step 3). But since accepting $C$ and $-C$ entails unacceptable conclusions ($\perp$), we either reject that (\ref{L}) shows \ref{liarL} is both true and false, and hence reject the dialetheic solution, or show the argument in (\ref{curry}) is invalid but the one in (\ref{L}) is not.

Since (\ref{curry}) only uses MP and CP if we want to reject it we have to reject MP or CP. MP is impeccable but if we reject it, since it is used in (\ref{L}), we also reject (\ref{L}).%, in which case (\ref{L}) does not show $L$ to be both true and false.%, which is what we aim to defend?.

Beall and Murzi \cite{beall} present a solution rejecting CP; we should reject that there is a deduction-theorem link between validity and conditionals, i.e., reject:
\begin{equation}\label{DT}\tag{DT}A\vdash B \textrm{ if and only if }A \rightarrow B\end{equation}
Since, if we accept \ref{DT} and MP, Curry's paradox cannot be solved, and if we reject \ref{DT} line 3 in (\ref{curry}) is not justified so the argument fails.

%OBJECTIONS: no justification for rejecting

However, rejecting \ref{DT} does not seem plausible; having $B$ as a consequence of $A$ is equivalent to having $A$ and $B$ joint by the consequence connective “$\rightarrow$". To make this more explicit, consider a sentence $\pi$ satisfying
%The argument from me to absurdity is valid.
$$ \pi \iff Val(\pi,\perp)$$
where $Val(A,B)$ means “The argument from $A$ to $B$ is valid", and an argument is valid if, when the premises are true, so is the conclusion. So, by definition 
\begin{equation}\label{V0}\tag{V0}Val(A,B) \rightarrow (A \rightarrow B).\end{equation}
%how can this hold if we just rejected DT!!!
Thus we have an analogous paradox:
\begin{equation}\label{vcurry}\end{equation}
1. $\pi \iff Val(\pi,\perp)$ (Given)\\
2.1 $\pi$ (Assumption for CP')\\
2.2 $Val(\pi,\perp)$ (1, 2.1; MP)\\
2.3 $\pi \rightarrow \perp$ (2.2; \ref{V0})\\
2.4 $\perp$ (2.1, 2.3; MP)\\
3. $Val(\pi,\perp)$ (2.1-2.4; CP')\\
4. $\pi\rightarrow \perp$ (3; \ref{V0})\\
5. $\pi$ (1, 3; MP)\\
6. $\perp$ (4, 5; MP)

This is \textit{v-Curry's paradox} \cite[p. 152]{beall}. Line 3 might seem to present the same problems given for line 3 in (\ref{curry}) but this one is justified since by showing $A\vdash B $ we are precisely showing the argument from $A$ to $B$ is valid; if $A\vdash B$ then $Val(A, B)$ (CP').

Since Curry's and v-Curry's paradoxes are essentially the same we should have a unified solution. But since for v-Curry \ref{V0} and CP' are justified we cannot reject (\ref{vcurry}) by rejecting those, so we cannot solve Curry's by rejecting CP. %to have a unified solution,

Therefore, since we cannot reject CP, MP is used in (\ref{L}), and since (\ref{curry}) only uses MP and CP, we cannot reject (\ref{curry}) at all, unless we reject (\ref{L}) too.

%So, if $A\vdash B$ then $Val(A, B)$ and if $Val(A,B)$ then  $(A \rightarrow B)$, hence the “if" part of \ref{DT} follows, which is the necessary part for running the conditonal proof in Curry's. Hence, rejecting CPs cannot been justified and the solution fails.

%Furthermore, since the two paradoxes are essentially the same we should expect to have a unified solution, hence if we rejected DT  the solution given by rejecting DT for both instances of the paradox. %is not plausible for either case and we reconsider the solution to the first one taking into account that we want to find a solution which applies to both instances.

Although it seems like (\ref{curry}) only uses MP and CP, Beall and Murzi \cite[p. 146]{beall} argue it implicitly uses \textit{Structural Contraction}:
\begin{equation}\label{SC}\tag{SC} \textrm{If }\Gamma\textrm{ }y,y \vdash b\textrm{ then }\Gamma\textrm{ } y \vdash b.\end{equation}
%That is, if $b$ follows from $y,y$ then it also follows from only $y$.
Other arguments for finding the truth-value of $C$ may use the \textit{rule of contraction}:
\begin{equation}\label{RC}\tag{RC}A \rightarrow (A \rightarrow B) \vdash A \rightarrow B \end{equation}
Thus, to reject the argument in (\ref{curry}) we can reject \ref{SC} and \ref{RC}. If we do, (\ref{curry}) is invalid since, by only one assumption of $C$ in 2.1 we discharged $C$ twice, in 2.2 and 2.3. Neither can we discharge $C$ only once in $(C \iff (C \rightarrow \perp)$ to conclude $\perp$ since we rejected \ref{RC}.%which is not valid if we reject \ref{SC}.

Similarly, in (\ref{vcurry}), only from one assumption of $\pi$ we discharge $\pi$ twice; in 2.3 and in 2.4. Hence, the argument is not valid if we reject \ref{SC}.

%reject the arguments in (\ref{curry}) and (\ref{vcurry}) and so
Therefore, by rejecting \ref{RC} and \ref{SC} we avoid Curry's and v-Curry's paradox in a unified way, and since neither \ref{RC} nor \ref{SC} is used in (\ref{L}) (except in steps 1.4 and 2.4 which can be omitted) we can still accept the dialetheic solution.

Rejecting \ref{SC} seems, however, very unintuitive and does not reflect how ordinary English is used. After all, if $A$ is true, everything following from $A$ is true and we can exploit $A$ as much as we like, it does not matter how many times we use $A$ in our arguments. 
%Furthermore, it is not clear what argument would make both $A$s in $A\rightarrow (A \rightarrow B)$ true. Surely, if one says, “If the earth is round then it is not flat and if the earth is round then it is not hyperbolic" (cambiar ejemplo) they do not need to prove twice that the earth is round for both consequents to follow. But, even if they did,  %Alternatively (the only option that seems to be left), it might be said that the sentence is not provable, which even less plausible.

Beall and Murzi \cite[p. 163]{beall} argue we are resistant to rejecting \ref{SC} because we have a view of validity as truth-preservation in all possible situations. But this might be abandoned if we think of premises instead of as partial descriptions of a situation, as resources. Then it is clearly different to have a single-$A$ resource from a double-$A$ resource. Accepting this requires changing our metaphysical account of validity.

Now, we clearly have $A \rightarrow A$ so, from one premise $A$, we can get as many $A$s as we like. But, of course, this is not true if we reject \ref{SC}; when we use our $A$ resource in $A \rightarrow A$ we obtain $A$ but we used one $A$ hence we are left with the same number of $A$ resources.

So, according to this new metaphysical account of validity, one ought to keep track of how many times they state a premise and how many they discharge. If someone had just proven the Riemann hypothesis and was presenting its implications they would have to state the Riemann hypothesis' truth, as well as any other premise, hundreds of times before presenting all the implications, which everyone, except maybe Beall and Murzi, would find pedantic. %In fact, they wouldn't get started with the proof before everyone topped listening since they would first have to prove the basics of algebra a few times in case they had ran out of those resources, making their proof invalid.

%And, even in doing so, it wouldn't be clear what makes having two Riemman hypothesis premise resources from having just one. If recreating the same argument counts as having doubled the resource then, since one can always recreate the same argument it always follows that if A then A and A. Hence, we can obtain a double-A resource from a single-A resource.

%corregir sc
Even if such a conception of validity is conceivable and we cannot from a single-A resource obtain a double-A resource, it can be argued it is possible to define some other concept which behaves like our natural understanding of validity and so satisfies \ref{SC}.%; an argument is valid' if, when the premises are true, the conclusion is true, and valid' follows .

Beall and Murzi consider this point for the case of connectives contraction. They argue conditionals do not contract as in \ref{RC} and no other connective does either. If there was one which did, say $+$, Curry's paradox would arise again for $a\iff(a+b)$. Thus, to avoid the paradox they reject such a connective exists. They could also argue that any concept satisfying \ref{SC} is meaningless.

To conclude, accepting the dialetheic solution forces us to accept $\perp$ unless we reject (\ref{curry}) by rejecting one of the rules of inference used. I ruled out rejecting MP and CP so the only possibility is to reject \ref{SC} and \ref{RC}, which forces us to change our metaphysical account of validity and accept that any paradox-free language is pedantic and has no contractible concepts. Even if this is an acceptable price to pay to solve the paradoxes, the Antinomy of the Liar and Curry's paradox have very similar structures so it should be possible to generalise their solutions to solve both and any of the same kind, since otherwise we risk the appearance of another paradox for which the specific solutions do not apply. Thus, the dialatheic solution is plausible but since it does not allow for a general solution we might want to reject it.% so if it is possible that other paradoxes will force we might want to reject it if other paradoxes arise. %argue like this: We should reject \ref{SC} hence we have to change our metaphysical account of validity; we should reject \ref{SC} hence we must reject any concept satisfying \ref{SC}; we should reject \ref{SC} hence English should be pedantic. However, maybe being pedantic is the price to pay for having a paradox-free language. 

\renewbibmacro*{in:}{%
  \iffieldundef{journaltitle}
    {}
    {\printtext{\bibstring{in}\intitlepunct}}}
    
\cleardoublepage
\renewcommand{\thepage}{}
\printbibliography

%Word count: 1492

\end{document}